\documentclass[reqno,b5paper]{amsart}
\usepackage{latexsym,graphicx,verbatim,color,cleveref,enumerate,fontenc,setspace,cite}
\usepackage{caption,amsthm,mathrsfs}
\usepackage{amsmath,amssymb}
\usepackage{subcaption}
\usepackage{enumerate}
\usepackage[mathscr]{eucal}
\newtheorem{thm}{Theorem}[section]
\newtheorem{deff}[thm]{Definition}
\newtheorem{defs}[thm]{Definitions}
\newtheorem{lemma}[thm]{Lemma}

\newtheorem{rem}[thm]{Remark}
\newtheorem{prop}[thm]{Proposition}
\newtheorem{cor}[thm]{Corollary}

\newtheorem {ex}[thm]{Example}

\newcommand{\erre}{\mbox{$\mathbb{R}$}}

\newcommand{\Fs}{\mbox{$\mathcal F$}}

\newcommand{\ttilde}{~}

\newcommand{\rFfreccia}{\mbox{$ \stackrel{\quad r_{\mathcal F} \quad}{\to}$}}
\newcommand{\oFfreccia}{\mbox{$ \stackrel{\quad o_{\mathcal F} \quad}{\to}$}}

\newcommand{\uno}{{\bf 1}}
\numberwithin{equation}{section}
\begin{document}
\begin{flushright}
%{\em \small Dedicated to Professor Paolo de Lucia \\  on the occasion of His eighty birthday}
\end{flushright}
\title[ $L^p$  spaces in vector lattices  ...]{\bf  $L^p$  spaces in vector lattices 
 and applications}
\author[A Boccuto \and D. Candeloro \and A. R. Sambucini]{\bf Antonio Boccuto \and 
Domenico Candeloro\and  Anna Rita Sambucini}
\date{March 16, 2016}
\newcommand{\acr}{\newline\indent}
\address{\llap{*\,}Department of Mathematics and Computer Sciences\acr
1, Via Vanvitelli  \acr 06123-I  Perugia \acr (Italy)}
\email{antonio.boccuto@unipg.it}
\address{\llap{*\,}Department of Mathematics and Computer Sciences\acr
1, Via Vanvitelli  \acr 06123-I  Perugia \acr (Italy)}
\email{domenico.candeloro@unipg.it}
\address{\llap{*\,}Department of Mathematics and Computer Sciences\acr
1, Via Vanvitelli  \acr 06123-I  Perugia \acr (Italy)}
\email{anna.sambucini@unipg.it}
\thanks{This work was supported by   University of Perugia - Department of Mathematics and 
Computer Sciences - Grant Nr 2010.011.0403 and by the Grant prot. U2014/000237 of GNAMPA - INDAM (Italy).}
\thanks{A. Boocuto orcid id: 0000-0003-3795-8856, D. Candeloro orcid id: 0000-0003-0526-5334, 
A.R. Sambucini orcid id: 0000-0003-0161-8729.}
\begin{abstract}
$L^p$ spaces are investigated for vector  lattice-valued functions, with respect to filter convergence.
As applications, some classical inequalities are extended to the vector lattice context, and some properties of
the Brownian Motion and the Brownian Bridge are studied, to solve some stochastic differential equations.
\end{abstract}
\maketitle
\section{Introduction}%\mg{bcs-16marzo2016} 
Function spaces, in particular $L^p$ spaces, play a central role in many 
problems in Mathematical Analysis and have lots of applications in several branches.  
The $L^p$ spaces are perhaps the most useful and important examples of Banach spaces; 
of independent and higher interest is the $L^2$ space,
whose origins are related with fundamental investigations and 
developments in Fourier analysis (\cite{Butzer, bm1, bm2}),  in  reconstruction  
of  signals,  integral  and  discrete  operators  (see for example \cite{AV1, AV2, 
BMV, BBV1, BBV2, BBV3, ccmv, CV, HOLMES, VZ, VZ1}) and  
Stochastic Integration (see also \cite{Bill, GL1, GL2, GLM, Grobler2, 
LETTA, s1990, s1990a, s1994, vw2014}).

In this paper some fundamental properties of $L^p$ spaces 
in the vector lattice setting are investigated, continuing a research
initiated by the authors in \cite{BC, taylor1, BC2011, bs1997, bs00,bcs-rl}
and developed later in \cite{dallas}.
The range of the involved functions is a vector lattice endowed with
filter/ideal convergence (for a related literature,
see also \cite{bd1, bd2, bbdm, bdlibro, bdpschurfilters, cs2014,dl2}).
Thanks to the triangle inequality, it is possible to view the space
$L^p$ as a metric space endowed with a distance of the type 
$d(f, g)= \|f - g\|_p$. Note that, in general, this space is not 
complete, as Example \ref{notcompleteness} shows.
As applications, several inequalities are given,
as well as some approximation results concerning processes related to the Brownian Motion.

The paper is organized as follows. In Section \ref{preliminaries}
basic notions and results are given, recalling 
some main properties of vector lattices, filters, modulars  and a Vitali-type Theorem.
In  Section \ref{seclp} some new results on $L^p$ spaces, for $p \in 
\mathbb{N}$ in the vector lattice context and a Minkowski inequality are given, together with an example 
in which it is shown that in general these spaces are not complete.
Section  \ref{app} is divided into two subsections; in the first 
one some inequalities  (Hermite-Hadamard, F\'{e}jer, Jensen and 
Schwartz) are proved. 
In the second subsection, applications to the  Brownian Motion and the Brownian Bridge are given, 
leading to the solution of some particular 
stochastic differential equations and to the reconstruction of a perturbed signal. 

%============================================
\section{Preliminaries}\label{preliminaries} %\mg{preliminaries}
 Some basic properties of vector lattices 
and filter convergence are recalled first. 
For the basic subjects and fundamental tools used in the vector lattice theory see, for instance,  \cite{LZ,Z}.
A vector lattice  $\mathbf{X}$ is said to be \textit{Dedekind complete} 
iff every nonempty subset $A \subset \mathbf{X}$, bounded from above, has a 
lattice supremum  in $\mathbf{X}$, denoted by $\bigvee A$.
From now on, $\mathbf{X}$ is a Dedekind complete vector lattice,  
$\mathbf{X}^+$ is the set of all strictly positive elements of
$\mathbf{X}$, and $\mathbf{X}^+_0=\mathbf{X}^+ \cup \{0\}$. For 
each $x \in \mathbf{X}$, let $|x| := x \vee (-x)$.
An extra element $+\infty$ will be added to $\mathbf{X}$, extending order 
and operations in a natural way, set
$\overline{\mathbf{X}}=\mathbf{X} \cup \{+ \infty\}$, 
$\overline{\mathbf{X}^+_0}=\mathbf{X}^+_0 \cup \{+ \infty\}$,
and assume $0 \cdot (+ \infty) = 0$.
A sequence $(p_n)_n$ in $\mathbf{X}$ is 
called \textit{$(o)$-sequence} 
iff it is decreasing and 
$\bigwedge_n p_n=0$.
An \textit{order unit} of $\mathbf{X}$ is an element $e$, such that for 
every $x \in \mathbf{X}$ there is 
a positive real number $c$ with $|x|\leq c e$. 

Let $\mathbf{X}_1$, $\mathbf{X}_2$, $\mathbf{X}$ be Dedekind 
complete vector lattices. We say that $(\mathbf{X}_1,\mathbf{X}_2,
\mathbf{X})$ is a \em product triple \rm 
iff a product $\cdot :\mathbf{X}_1
\times \mathbf{X}_2 \to \mathbf{X}$ is defined, satisfying natural 
conditions of compatibility  (see for instance \cite[Assumption 2.1]{BC2011}). 

\begin{rem}\label{mov} %\mg{mov} 
\rm
A vector lattice  $\mathbf{X}$ is called an \textit{$f$-algebra} 
(see also \cite[Definition 140.8]{Z}) iff 
there exists in $\mathbf{X}$ an associative multiplication, satisfying 
the usual algebraic 
properties, with  $xy \geq 0$ whenever $x \geq 0$ and $y \geq 0$
and such that $x \wedge y 
= 0$ implies $(x \cdot z) \wedge y = 0$ whenever 
$x$, $y \in \mathbf{X}$ and $z 
\in \mathbf{X}_0^+$.
%In this paper we always assume that
The following condition will be required in the paper.
\begin{itemize}
\item[ \textbf{($H_0$)}]
 $(\mathbf{X}_1, \mathbf{X}_2, \mathbf{X})$ is a product triple,  
and $\mathbf{X}, \mathbf{X}_1$ are endowed with order units 
$e$, $e_1$ respectively.
\end{itemize}
Note that every lattice $\mathbf{X}$ 
equipped with an order unit is an $f$-algebra. 
Indeed, by the Maeda-Ogasawara-Vulikh representation theorem
\ref{movth} (see also \cite{WRIGHT}), 
$\mathbf{X}$ is algebraically and lattice isomorphic to 
the space $C(\Omega)$ of all continuous real-valued functions
defined on a suitable compact and extremely disconnected 
topological space $\Omega$. 
So, it is not difficult to deduce that $\mathbf{X}$ is an
$f$-algebra since $\mathbb{R}$ is.
\end{rem}

Given any fixed countable set $Z$, a class  $\mathcal{F}$ of subsets of $Z$
is called  a \textit{filter} of $Z$ iff $\emptyset \not \in {\mathcal F}$,
$A \cap B \in {\mathcal F}$ whenever $A$, $B \in \mathcal{F}$ and for each
$A \in \mathcal{F}$ and $B \supset A$ it is $B\in \mathcal{F}$.
The symbol  ${\mathcal F}_{\text{cofin}}$ denotes  the filter of all cofinite subsets of $Z$.
A filter of $Z$ is said to be \textit{free} iff it contains ${\mathcal F}_{\text{cofin}}$.
An example of free filter is the filter
of all subsets of $\mathbb{N}$,
having asymptotic density one.
If ${\mathcal F}$ is a filter of $Z$, then let  ${\mathcal F}
\otimes {\mathcal F}$  be the \textit{product filter}
of $Z \times Z$, defined by
% \mg{productfilter} 
\begin{align}\label{productfilter}
{\mathcal F} \otimes {\mathcal F} :=\{C
\subset Z \times Z: \text{  there exist  } A, B \in 
{\mathcal F} \text{  with  }
A \times B \subset C \}. 
\end{align}
\begin{deff}\label{filterconv} \rm %\mg{filterconv}
Let ${\mathcal F}$ be any filter of $Z$.
A sequence $(x_z)_{z \in Z}$ in $\mathbf{X}$ \textit{$(o_{\mathcal F})$-converges to $x
\in R$}    ($x_z \oFfreccia x$) iff there exists an $(o)$-sequence $(\sigma_p)_p $ in $\mathbf{X}$  such
that for all $p \in \mathbb{N}$ the set $\{ z \in Z: |x_z - x| \leq \sigma_p  \}$ belongs to
${\mathcal F}$.\acr
A sequence $(x_z)_{z \in Z}$ in $\mathbf{X}$ \textit{$(r_{\mathcal F})$-converges to $x
\in \mathbf{X}$}  ($x_z \rFfreccia x$) iff there exist a 
$u \in \mathbf{X}^+$ and an $(o)$-sequence $(\varepsilon_p)_p$ in 
$\mathbb{R}^+$  such
that for every 
$p \in \mathbb{N}$ the set 
$\{ z \in Z: |x_z - x| \leq  \varepsilon_p  u\}$
is an element of ${\mathcal F}$.\acr
A sequence $(x_z)_z$ in $\mathbf{X}$ \textit{$(o)$}-converges to $x \in \mathbf{X}$ 
(in the classical sense) iff
it $(o_{{\mathcal F}_{\text{cofin}}})$-converges to $x$
(see also \cite{bdpschurfilters}).\\
\end{deff}

Let $G$ be any infinite set, ${\mathcal P}(G)$ be the family of all
subsets of $G$, ${\mathcal A}
\subset {\mathcal P}(G)$ be an algebra and $\mu:{\mathcal A}
\to \overline{(\mathbf{X}_2)^+_0}$ be a finitely additive measure.
The symbol  $\mu^*$ denotes 
the \textit{outer measure} associated to $\mu$, namely
$\displaystyle{\mu^*(B):=
\wedge_{A\in {\mathcal A},A\supset B} \mu(A)}$, 
$B \in {\mathcal P}(G)$, and 
${\mathcal A_b}$ is the family of the sets $B\in 
{\mathcal A}$ with $\mu(B)\in \mathbf{X}_2$.
For every $A \in {\mathcal A}$ and $u \in \mathbf{X}_1$, 
let $u\cdot 1_A : G \to \mathbf{X}_1$ be the function whose 
values are $u$ when $t\in A$ and 0 otherwise. 

\noindent
As in \cite[Subsection 2.1]{dallas}, 
the modulars in the 
vector lattice setting are introduced
(for the classical case and a related literature, see e.g. 
\cite{BMV, KOZLOWSKI,Musielak,RAO}). \acr
Let $T$ be a linear sublattice of $\mathbf{X}_1^G$,
such that $e_1 \cdot 1_A \in T$ for every $A \in {\mathcal A_b}$.\acr

\noindent
A functional $\rho:T \to \overline{\mathbf{X}^+_0}$ is said to be a
\textit{modular} on $T$ iff it satisfies the following properties.
\begin{itemize}
\item[($\rho_0$)] $\rho(0)=0$;
\item[($\rho_1$)] $\rho(-f)=\rho(f)$ for every $f \in T$;
\item[($\rho_{2}$)] $\rho(\alpha_1 f + \alpha_2 h) \leq \rho(f) + 
\rho(h)$ for every
$f$, $h \in T$ and for any $\alpha_1$, $\alpha_2 \geq 0$  with
$\alpha_1 + \alpha_2 =1$.
\end{itemize}
Moreover, some additional conditions will be required.
\begin{itemize}
\item[($\rho_m$)] A modular $\rho$ is \textit{monotone} iff $\rho(f)
\leq \rho(h)$ for every $f$, $h \in T$ with
$|f|\leq |h|$. In this case, if $f \in T$, then $|f| \in T$ and hence $\rho(f) = \rho(|f|)$.
\item[($\rho_{co}$) ]
A  modular $\rho$ is \textit{convex} iff
$\rho(\alpha_1 f_1 + \alpha_2 f_2) \leq \alpha_1 \rho(f_1) +
\alpha_2 \rho(f_2)$ for all
$f_1$, $f_2 \in T$ and for any real numbers $\alpha_1$, $\alpha_2 \geq 0$ with
$\alpha_1 + \alpha_2=1$.
\item[($\rho_f$)]
A modular $\rho$ is \textit{finite} iff 
for every $A \in {\mathcal A_b}$  and every $(o)$-sequence $(\varepsilon_p)_p$ in $\mathbb{R}^+$,
 the sequence $(\rho(e_1 \varepsilon_p  1_A))_p$ is $(r)$-convergent to $0$
(for the case $\mathbf{X}=\mathbf{X}_1=\mathbf{X}_2=
\mathbb{R}$, see for instance \cite{BMV}).
\end{itemize}
 The concept of (equi-) absolute continuity in the context of 
modulars and filter convergence is introduced here. 
\begin{itemize}
\item[($a_{\rho}$)]
A map $f \in T$ is \textit{$(o_{\mathcal F})$-absolutely
continuous with respect to the modular $\rho$} 
(shortly, \textit{absolutely continuous}) iff there is a positive real constant $\alpha$,
satisfying the following properties.
\begin{itemize}
\item[($a_{\rho}(1)$)] For each $(o)$-sequence $(\sigma_p)_p$ in $ 
\mathbf{X}_2^+$ there exists an $(o)$-sequence $(w_p)_p$ in $R^+$
such that for all $p \in \mathbb{N}$ and whenever $\mu(B) \leq \sigma_p$ it is
$\rho(\alpha f  1_B) \leq w_p$;
\item[($a_{\rho}(2)$)] there is an $(o)$-sequence $(z_m)_m$ in $\mathbf{X}^+$ such that to
each $m \in \mathbb{N}$ there corresponds a set $B_m \in {\mathcal A_b}$
with $\rho(\alpha f  1_{G \setminus B_m}) \leq z_m$.
\end{itemize}
\item[($ac_{\rho}$)]
Given a modular $\rho$ and any free filter ${\mathcal F}$ of $Z$,
 a sequence $f_z:G \to \mathbb{R}$, $z \in Z$, is said to be
\textit{$\rho$-${\mathcal F}$-equi-absolutely continuous}, or in short
\textit{equi-absolutely continuous,} iff
there is  $\alpha \in \mathbb{R}^+$, satisfying the following two conditions.
\begin{itemize}
\item[($ac_{\rho}(1)$)] For every $(o)$-sequence $(\sigma_p)_p$ in 
$\mathbf{X}_2^+$  there are an $(o)$-sequence $(w_p)_p$ in $\mathbf{X}^+$
and a sequence $(\Lambda^{p})_p \in {\mathcal F}$ with
$\rho(\alpha f_z  1_B) \leq w_p$
whenever $z \in \Lambda^{p}$ and $\mu(B) \leq \sigma_p$, $p\in \mathbb{N}$;
\item[($ac_{\rho}(2)$)] there are an $(o)$-sequence $(r_m)_m$ in $\mathbf{X}^+$ 
and a sequence $(B_m)_m \in {\mathcal A_b}$
such that,  for all $m\in \mathbb{N}$, it is
%\mg{equiac2} 
\begin{align}\label{equiac2}
\Lambda^{m}:=\{z \in Z: \rho(\alpha f_z  1_{G \setminus B_m})
\leq r_m \} \in {\mathcal F}.
\end{align}
\end{itemize}
\end{itemize}
The concepts of filter uniform 
convergence and convergence in measure are recalled  (see also  \cite{BC}).
Let ${\mathcal F}$ be any fixed free filter of $Z$.
\begin{defs} \rm \label{u-misura} %\mg{u-misura}
\ttilde
\begin{itemize} \vskip-1cm
\item[\ref{u-misura}.1)]  A sequence of functions 
$(f_z)_{z \in Z}$ in 
${\mathbf{X}_1}^G$ is said to
{\em  \textit{ $(r_{\mathcal F})$}-converge uniformly} (shortly, 
\textit{converge uniformly}) to $f$,
iff there exists an $(o)$-sequence $(\varepsilon_p)_p$ in $\mathbb{R}^+$ with
\[ \{ z \in Z : \bigvee_{t \in G} \, |f_z(t)-f(t)| 
\leq \varepsilon_p \, e_1 \} \in {\mathcal F}
\quad \text{for  any  } p \in \mathbb{N}.\] In this case one can write:
$$(r_{\mathcal F})\lim_z \Bigl( \bigvee_{t \in G}  |f_z(t)-f(t)|\Bigr)=0.$$
\item[\ref{u-misura}.2)] Given a sequence  $(f_z)_z$ in 
${\mathbf{X}}_1^G$ and $f \in {\mathbf{X}}_1^G$, 
we say that
$(f_z)_z$ \textit{$(r_{\mathcal F})$-converges in measure}
(shortly, \textit{$\mu$-converges}) to $f$,
iff there are  two $(o)$-sequences,
$(\varepsilon_p)_p$ in $\mathbb{R}^+$, 
$(\sigma_p)_p$  in ${\mathbf{X}}_2^+$,
and a double sequence $(A_z^p)_{(z,p) \in Z 
\times \mathbb{N}}$ in ${\mathcal A}$ such that
$A_z^p \supset \{t \in G: |f_z(t)-f(t)| \not\leq \varepsilon_p e_1\}$ for every $z \in Z$ 
and $p \in \mathbb{N}$, and 
$\{z \in Z:  \mu(A_z^p) \leq \sigma_p \} \in {\mathcal F}$  for all $p \in \mathbb{N}$. 
\end{itemize}
\end{defs}
Using these assumptions and notations, in \cite[Theorem 2.3]{dallas}
a Vitali-type theorem was obtained. For a historical overview on this topic see also
\cite{bdlibro,CHOKSI} and their bibliographies.
%\mg{vitali} 
\begin{thm}\label{vitali}  
{\rm(Vitali)}
Let $\rho$ be a monotone and  finite modular, and 
$\mathcal F$ be a fixed free filter of $Z$. If 
$(f_z)_z$ is a sequence in $T$, 
$\mu$-convergent  to $0$ and  equi-absolutely continuous, then there 
is a positive real number $\alpha$ with
$$(o_{\mathcal F})\lim_z \rho(\alpha  f_z)=0.$$
\end{thm}

\noindent As an application,  a Cauchy-type property for 
$\rho$-convergence of function sequences can be obtained. 
\begin{thm}\label{lebesguebis}  %\mg{lebesguebis}
Let ${\mathcal F}$ and $\rho$  be as in Theorem 
{\rm \ref{vitali}}, and suppose
that $(f_n)_{n\in \mathbb{N}}$ 
is a sequence in $T$, such that the double sequence
$(f_h-f_q)_{h,q}$
$(o_{{\mathcal F}\otimes{\mathcal F}})$-converges 
in measure to $0$. If 
$(f_n)_n$ is equi-absolutely continuous, then
there is $\alpha > 0$ with $\displaystyle{
(o_{{\mathcal F} \otimes {\mathcal F}})\lim_{h,q} \rho(\alpha  (f_h-
f_q))=0}.$
\end{thm}
\begin{proof} It is enough to
replace, in Theorem \ref{vitali}, $Z$ with $\mathbb{N}
\times \mathbb{N}$ and 
${\mathcal F}$ with ${\mathcal F}\otimes {\mathcal F}$,
respectively.
\end{proof}

As a consequence it follows
\begin{cor} {\label{cauchydominato}} %\mg{cauchydominato}
Let ${\mathcal F}$ and $\rho$ be as in 
Theorem {\rm \ref{vitali}}, and suppose 
that $(f_n)_n$ is 
$(o_{\mathcal F})$-convergent in measure to 0.
If there exist an absolutely continuous
function $g$ in $T$ and an element $F_0 \in {\mathcal F}$, 
such that
$|f_z(t)| \leq g(t)$ for all $z \in F_0$ and $t \in G$,
then there is a positive real number $\alpha$ with 
$(o_{\mathcal F})\lim_z \rho(\alpha \, f_z)=0.$
\end{cor}
%============= inizio applicazioni=============================
\section{$L^p$ spaces}\label{seclp} %\mg{seclp}
Some definitions 
which will be used in the sequel are recalled here for the sake of simplicity.
\begin{deff}\label{d3.2}\rm \cite[Definition 3.2]{dallas} %\mg{d3.2}
A function $f \in \mathbf{X}_1^G$ 
is said to  be \textit{simple} iff $f(G)$ is a finite set and
$f^{-1} (\{x\}) \in {\mathcal A}$ 
for every $x \in \mathbf{X}_1$. The space of all simple functions 
is denoted by ${\mathcal S}$.
\end{deff}

Let $L^*$ be the set of all simple functions $f\in {\mathcal S}$ vanishing outside a set of 
finite $\mu$-measure.
If $f \in L^*$,  its usual integral is denoted by 
$\displaystyle{\int_G f(t)  d\mu(t)}$.
It is not difficult to see 
that the functional $\iota: L^* \to \mathbf{X}$ defined as
%\mg{modulareintegrale} 
\begin{eqnarray}\label{modulareintegrale}
\iota(f) :=\int_G |f(t)|  d\mu(t), \quad f \in L^*,
\end{eqnarray}
is a monotone finite modular, and it is also linear and  additive on 
positive functions and constants 
(see also \cite[Remark 3.3]{dallas}).
\begin{deff}\label{integrabilitainfinita} \rm (\cite[Definition 3.4]{dallas}) %\mg{integrabilitainfinita}
A positive function $f \in {\mathbf{X}}_1^G$ is
\textit{integrable} iff there exist an equi-absolutely continuous sequence
of functions $(f_n)_n$ in $L^*$, $\mu$-convergent to $f$, and a map 
$l:{\mathcal A} \to \mathbf{X}$, with
%\mg{zvezda} 
\begin{eqnarray}\label{zvezda} (o_{\mathcal  F})\lim_n 
\bigvee_{A \in {\mathcal A}}
\Bigl|\int_A f_n(t)  d\mu(t)-l(A)\Bigr|=0
\end{eqnarray}
(the sequence $(f_n)_n$ is said to be \textit{defining}). In this case, we say that
$$l(A)=(o_{\mathcal F})\lim_n \int_A    f_n(t) d\mu(t)$$ 
\textit{uniformly with respect to} $ A \in {\mathcal A}$. 
\acr
Note that $l(A)$ is independent of the choice of the defining sequence.
If $f \in {\mathbf{X}}_1^G$,  then  $f$ is said to be
\textit{integrable} iff the functions $t \mapsto f(t) \vee 0$ and
$t \mapsto (-f(t)) \vee 0$ are integrable.\acr
\end{deff}

The $L^p$ spaces in the vector lattice setting will be introduced now.
Assume that 
$\mathbf{X}_1 = \mathbf{X}$,  ( $\mathbf{X}$ is a lattice ordered algebra with
complete multiplication, see also \cite{Z}), 
$\mathbf{X}_2= \mathbb{R}$.

\begin{deff}\label{lp} %\mg{lp}
\rm
Let $f \in \mathbf{X}^G $.
We say that $f \in L^p$ iff both $f$ and 
$f^p$ belong to $L$ 
according with Definition \ref{integrabilitainfinita} 
with a common basic defining sequence $(f_n)_n$ (this means that 
$(f_n)_n$ is a defining sequence for $f$ and $(f_n^p)_n$ is a defining 
sequence for $f^p$, respectively).
\end{deff}
It is not difficult to see that $L^p$ is a linear space. Since homogeneity 
is straightforward, we just prove that $L^p$ is stable with respect to 
addition. The  $\sigma$-finiteness 
property for measurable functions will be proved first.
\begin{lemma}\label{l1} %\mg{l1}
Let $(f_n)_n$ be a sequence of simple functions, $\mu$-converging to a mapping $f$. 
Then there exist an $(o)$-sequence $(\beta_k)_k$ in 
$\mathbb{R}^+$, an increasing sequence $(N_k)_k$ of positive integers and a sequence $(H_k)_k$ in ${\mathcal A}$ with  
 $\mu(H_k) \leq \beta_k$ and $\{ t \in G: |f(t)| \not\leq N_k e \} \subset H_k$ for every $k\in\mathbb{N}$.
\end{lemma}
\begin{proof}
Thanks to $\mu$-convergence, there exist an $(o)$-sequence $(\sigma_k)_k$ in 
$\mathbb{R}^+$, an $(o)$-sequence $(\varepsilon_k)_k$ in the real interval $]0,1]$ and a sequence 
$(F^k)_k$ in $\mathcal{F}$ such that, for every integer $k$ and every $z \in F^k$ there exists an element 
$A_z^k \in {\mathcal A}$ 
satisfying $\mu(A_z^k) \leq \sigma_p$ and $\{t \in G: |f(t) - f_{z} (t)| \not\leq \varepsilon_k e \} \subset A_z^k$. \acr
Now for every $k$ let us denote by $z_k$ the least element of $F^k$ 
and by 
$N_k$ any positive integer such that $|f_{z_k} (t)| \leq N_k e$ for every 
$t \in G$, and set $H_k:=A^k_{z_k}$. Then for every integer $k$ 
it is $\mu(H_k) \leq \sigma_k$ and 
$\{ t \in G: |f(t) - f_{z_k} (t) | \not\leq e\} \subset H_k,$ and therefore
$\{ t \in G: |f(t)| \not\leq (N_k +1) e \} \subset H_k,$
from which the assertion follows, just replacing $\sigma_k$ 
with $\beta_k$ 
and $N_k +1$ with $N_k$.
\end{proof}
\begin{prop}\label{sum-p} %\mg{sum-p}
If $f$, $g \in L^p$, then $f+g \in L^p$.
\end{prop}
\begin{proof}
Let $(f_n)_n$ and $(g_n)_n$ be two defining sequences for $f$, $g$, respectively. 
From the properties of $\mu$-convergence it follows that the sequence $(f_n + g_n)_n$ is $\mu$-convergent to $f+g$. 
Moreover, since $(f_n)_n$ and $(g_n)_n$ are equi-absolutely continuous 
(with respect to the modular $\iota$), 
it is not difficult to check that $(f_n + g_n)_n$ is too. 
\acr
Indeed, if $(\sigma_k)_k$ is any $(o)$-sequence in 
%$\mathbf{X}^+$: 
$\mathbb{R}^+$, then there are two $(o)$-sequences $(w_k)_k$ and 
$(w_k^{'})_k$ in $\mathbf{X}$ and two sequences $(\Xi_k)_k, 
(\Xi_k^{'})_k$ in ${\mathcal F}$, such that for every 
$k \in \mathbb{N}$ it is
$\iota(f_z 1_B) \leq w_k$,  $\iota(g_k 1_B) \leq w_k^{'}$
as soon as $z \in \Xi_k \cap \Xi_k^{'}$ and $\mu(B) \leq \sigma_k$. 
Thus, choosing $w_k^{*} = w_k + w_k^{'}$ and
$\Xi_k^{*} =\Xi_k \cap \Xi_k^{'}$, it is
$\iota((f_n + g_n)1_B) \leq w_k^{*} $
as soon as $\mu(B) \leq \sigma_k$ and $z \in \Xi_k^{*}$. This proves property ($ac_{\rho}(1)$). 
\acr Moreover, there are two $(o)$-sequences 
$(r_m)_m, (r_m^{'})_m$ in $\mathbf{X}^+$ and two sequences $(B_m)_m, (B_m^{'})_m$ 
in ${\mathcal A}_b$ such that the sets
$\Lambda_m := \{ z \in Z : \iota(f_z 1_{G \setminus B_m} \leq r_m \},
\quad \Lambda_m^{'} :=
 \{ z \in Z : \iota(g_z  1_{G \setminus B_m^{'}} \leq r_m^{'} \}$
belong to ${\mathcal F}$. So, 
taking  $r_m^{*} = r_m + r_m^{'}$ and
$B_m^{*} =B_m \cup B_m^{'}$, for every $z \in \Lambda_m \cap \Lambda_m^{'}$ it is
\[\iota((f_z + g_z)  1_{G \setminus B_m^{*}})  \leq 
\iota(f_z  1_{G \setminus B_m^{*}}) + \iota(g_z  1_{G \setminus B_m^{*}}) \leq
\iota(f_z  1_{G \setminus B_m}) + \iota(g_z  1_{G \setminus B_m^{'}}) \leq r_m^{*},\]
which proves ($ac_{\rho}(2)$). \acr
Now for the same reason, since $(f_n^p)_n$ and $(g_n^p)_n$ are 
defining sequences for $f^p$, $g^p$, respectively, we deduce that $(|
f_n|^p + |g_n|^p)_n$ is an equi-absolutely continuous sequence (we 
recall that the absolute continuity is essentially a condition on $|f|$). 
From this, since
$ |f_n + g_n|^p \leq (|f_n| + |g_n|)^p \leq 2^p(|f_n| \vee |g_n|)^p 
\leq 2^p (|f_n|^p + |g_n|^p),
$
 it is clear that the sequence $(f_n + g_n)^p$ is equi-absolutely 
continuous. The next step is to prove that the sequence $((f_n + 
g_n)^p)_n$ is $\mu$-convergent to $(f+g)^p$. 
Thanks to $\mu$-convergence of the four sequences $(f_n)_n$, 
$(g_n)$, $(f_n^p)_n$, $(g_n^p)_n$, there exist an $(o)$-sequence 
$(\varepsilon_k)_k$ in $]0,1[$ and a sequence
$(F_k)_k$ in ${\mathcal F}$ 
such that for every integer $k$ and every $z \in F_k$ there is $A_k^z \in 
{\mathcal A}$ such that
$\mu(A_k^z) \leq \sigma_p$, and furthermore 
\begin{eqnarray*}
&& \{ t \in G: |f_z(t) - f(t)| \not\leq \varepsilon_k e \} \cup
 \{ t \in G: |g_z(t) - g(t)| \not\leq \varepsilon_k e \} \subset A_k^z,\\
&&
\{ t \in G: |f^p_z(t) - f^p(t)| \not\leq \varepsilon_k e \} \cup
 \{ t \in G: |g^p_z(t) - g^p(t)| \not\leq \varepsilon_k e \} \subset A_k^z.
\end{eqnarray*}
As in Lemma \ref{l1}, without loss of generality, the quantities $\beta_k$, $N_k$, $H_k$
 here obtained can be considered the same for both $f,g$.\acr
Now  a subsequence $ (\varepsilon^{'}_k)_k$  of $ (\varepsilon_k)_k$ will be found such that
$\varepsilon^{'}_k \leq (k p 2^p (N_k +1)^p)^{-1}$
for every $k$. In correspondence with $(\varepsilon^{'}_k)_k$, we denote by
$(\sigma^{'}_k)_k$ and $(F^{'}_k)_k$ the subsequences of 
$(\sigma_k)_k$ and $(F_k)_k$, respectively. Fix $k \in \mathbb{N}$, 
choose an element $z \in F^{'}_k$ and let $t \not\in A_k^z \cup H_k$. 
Then it is
\begin{eqnarray*}
&& (f_z (t) + g_z (t))^p - (f (t) + g(t))^p = [ (f_z(t) - f(t)) +(g_z(t) - g(t)] \cdot\\  &&
\{  (f_z (t) + g_z (t))^{p-1} +  (f_z (t) + g_z (t))^{p-2} (f(t) + g(t)) + 
\ldots + (f(t) + g(t))^{p-1} \},
\end{eqnarray*}
and therefore 
$ |(f_z (t) + g_z (t))^p - (f (t) + g(t))^p| \leq 2 \varepsilon^{'}_k p 2^p (N_k+1)^p \leq \dfrac{2}{k} e \downarrow 0.
$
Since $\mu(A_k^z \cup H_k) \leq \beta_k + \sigma^{'}_k$, this is 
sufficient to conclude the proof of $\mu$-convergence. Thus also the 
double sequence $(n,k) \mapsto (f_n +g_n)^p - (f_k + g_k)^p$ is 
equi-absolutely continuous and $\mu$-convergent to 0 (with respect to 
the product filter $\Fs \otimes \Fs$).
From this, thanks to \cite[Theorem 2.3]{dallas}, 
the integrals
$$\int_G |(f_n(t) +g_n(t))^p - (f_k(t) + g_k(t))^p| d\mu(t)$$ 
converge to $0$ as $n,k \to \infty$, and so 
$$\int_E (f_n(t) +g_n(t))^p - (f_k(t) + g_k(t))^p d\mu(t)$$ 
converge to 0 uniformly with respect to $E$, which in turn implies that the limit
\[(o_{\mathcal F})\lim_n \int_E (f_n(t) + g_n(t))^p d\mu(t)\]
exists uniformly with respect to $E$ (see 
also \cite[Proposition 2.14]{bd2}), thus proving that $((f_n + 
g_n)^p)_n$ is a defining sequence for $(f+g)^p$.
\end{proof}

One easily sees that, if $f \in L^p$, also $|f|$ is and the map 
$\displaystyle{f \mapsto \int_G |f|^p d\mu}$ is 
a monotone and finite modular in $L^p$. \acr

A norm in the space $L^p$ can be defined in the following way:
\[ \| f\|_p := \inf \Bigl\{ \varepsilon > 0 : \int_G \Bigl(\dfrac{|f(t)|}{\varepsilon} \Bigr)^p d\mu(t) \leq e \Bigr\}
=\Bigl\| \int_G |f(t)|^p d\mu(t)\Bigr\|_e^{\frac{1}{p}},
\]
where $\|x\|_e = \inf \{ \varepsilon > 0: |x| \leq \varepsilon e \}$ is 
the $M$-norm in $\mathbf{X}$ associated with $e$
(see also \cite{SISY}).\acr

It is not difficult to verify   that $\| \cdot \|_p$ is positively 
homogeneous; in order to prove the Minkowski 
inequality  the  Maeda-Ogasawara-Vulikh representation 
theorem for vector lattices is needed (see also \cite{WF, LZ, WRIGHT}).
\begin{thm} \label{movth} {\rm (MOV
representation \cite[\S1]{WRIGHT})}
Given a Dedekind complete vector lattice $\mathbf{X}$ 
with order unit $e$, there
exists a compact extremely disconnected topological space
$\Omega$, unique up to homeomorphisms, such that $\mathbf{X}$ 
is algebraically and lattice isomorphic to 
$C(\Omega):=\{h \in {\mathbb{R}}^{\Omega}:h$
is continuous$\}$. 
\end{thm}
 
This allows to prove the following
\begin{prop}\label{mink} %\mg{mink}
{\rm (Minkowski inequality)}
 For every $p \in \mathbb{N}$  and $f,g \in L^p$, it is
\[ \| f+ g\|_p \leq \|f\|_p + \|g\|_p.\]
\end{prop}
\begin{proof}
Let $\Omega$ be as in Theorem \ref{movth}.
By Remark \ref{mov} and \cite[Lemma]{mali},  
for every $\alpha \in ]0,1[$, $t\in G$ and
$\omega \in \Omega$ 
it is
%\mg{vom} 
\begin{eqnarray}\label{vom}
(|f(t)(\omega)|+|g(t)(\omega)|)^p \leq 
\alpha^{1-p}|f(t)(\omega)|^p + (1-\alpha)^{1-p} 
|g(t)(\omega)|^p.
\end{eqnarray}
Since, 
by \cite[Corollary 2.20]{WF} and Theorem \ref{movth} the 
representation preserves the multiplication, 
from  (\ref{vom}) it is
for each $\alpha \in ]0,1[$ and $t\in G$ it is
\[ (|f(t)|+|g(t)|)^p \leq \alpha^{1-p}|f(t)|^p + (1-\alpha)^{1-p} 
|g(t)|^p.\]
So, following again \cite{mali}, it is
\begin{eqnarray*}
\int_G | f(t)+g(t)|^p d\mu(t) &\leq& \int_G (| f(t)|+|g(t)|)^p d\mu(t) \leq \\ &\leq&
\int_G (\alpha^{1-p} |f(t)|^p + (1-\alpha)^{1-p} |g(t)|^p ) d\mu(t) = \\ &=&
\alpha^{1-p} \int_G  |f(t)|^p  d\mu(t) + (1-\alpha)^{1-p}  \int_G  |g(t)|^p d\mu(t) \leq \\ &\leq&
 \int_G  |f(t)|^p  d\mu(t) +  \int_G  |g(t)|^p d\mu(t)
\end{eqnarray*}
Now, since 
\[ \| f + g \|_p =\Bigl\| \int_G |f(t)+ g(t)|^p d\mu(t)\Bigr\|_e^{\frac{1}{p}} \]
and by the monotonicity of the $M$-norm, the conclusion follows.
\end{proof}

Moreover one can see that, when
$\|f\|_p = 0$, then $\displaystyle{\int_G |f(t)|^p d\mu(t)=0}$, 
which implies that $\mu(E) 
=0$ for each $\varepsilon > 0$ and every $E \in {\mathcal A}$ such that
$f(t) \geq \varepsilon e$ for all $t \in E$.
This property can be expressed by saying that $f$ is \textit{essentially $\mu$-null}.
The essentially $\mu$-null functions in $L^p$ form a subspace, 
which allows us to introduce an equivalence relation in $L^p$,
by setting
\[ f \sim g \quad \mbox{iff} \quad f-g  \quad \mbox{ is essentially} \quad \mu\mbox{-null}.\]
 There is a large literature on completeness of $L^p$ spaces, see for example \cite{dl1,bb2000}. In this  setting
the space $L^p$ is not complete in general, as the following example shows.
%\mg{da controllare}
\begin{ex}\label{notcompleteness} \rm 
Following \cite[Example 2.2]{bb2000},
let $G = \mathbb{N}$ and $\mathcal{A}$ be the family of all finite 
or cofinite sets.
Observe that the completion of $\mathcal{A}$ in the sense of 
\cite[Section 1.6]{bb2000} is $\mathcal{P}(\mathbb{N})$.
Let
%\mg{mumu}
\begin{eqnarray}\label{mumu} 
\mu(A) := \left\{ \begin{array}{ll}
\mu(A) = \displaystyle{\sum_{i \in A} 2^{-i}}, &  \text{if } 
A \text{ is finite},\\ & \\
\mu(A) = 2 - \displaystyle{\sum_{i \in A^c} 2^{-i}}, 
& \text{if } A \text{ is cofinite}.
\end{array} \right. 
\end{eqnarray}
In this setting, the $\mu$-convergence of a sequence $f_n$ to $f$ can 
be expressed in the following way:
for every $\varepsilon > 0$ there exists $n(\varepsilon) \in \mathbb{N}$ such that, for every 
$n > n(\varepsilon)$, it is $\mu^*(\{s \in \mathbb{N} : |f_n(s) -f(s)| \not\leq \varepsilon e \}) \leq \varepsilon$,
where $\mu^*$ denotes the usual outer measure.\\
 Let $A_n = \{1,2, \ldots, n\}$, $n\in\mathbb{N}$. 
The sequence $(\uno_{A_n})_n \subset  L^*$ and then $\uno_{A_n} \in L^p$ for every $p \in \mathbb{N}$. 
Note that  $(\uno_{A_n})_n$ is a Cauchy sequence in $L^p$. In fact, for every $m,n \in \mathbb{N}$, it is
$\int_{\mathbb{N}} | \uno_{A_n}(s) - 
\uno_{A_m}(s)|^p d\mu(s) = e^p \mu(A_n \triangle A_m).$
Suppose, by contradiction, that $L^p$ is complete, so there exists $f  \in L^1$ such that $(\uno_{A_n})_n$ 
converges to $f$ in $L^1$.
Following \cite[Proposition 1.8 (b) and Lemma 2.1]{bb2000} we will show that $f = \uno_A$ for some $A \subset \mathbb{N}$.
Let $(n_k)_k$ be an increasing sequence in $\mathbb{N}$ such that
$\mu^* (\{ s \in \mathbb{N}: |\uno_{A_{n_k}} (s) - f(s)| \not\leq k^{-1} \}) \leq  k^{-1} e.$
Let $B_k = \{ s \in \mathbb{N}: |\uno_{A_{n_k}}(s) - f(s)| \not\leq k^{-1} e \}$.
As in \cite[Lemma 2.1]{bb2000}, it follows that
$ s \in A_{n_k} \cap B_k^c \Longrightarrow | \uno_{\mathbb{N}}(s) - f(s)| \leq  k^{-1} e; \quad
 s \in A^c_{n_k} \cap B_k^c \Longrightarrow |f(s)| \leq k^{-1} e.$
Let now $A:= \displaystyle{ \cup_{k=3}^{\infty}}
(A_{n_k}  \cap B_k^c)$,  
then $A \triangle A_{n_k} \subset B_k $, and hence
$\int_{\mathbb{N}} | \uno_{A_n}(s)-\uno_A(s)| d\mu(s) = \mu(A_n 
\triangle A) e \to 0.$
So, $f \sim \uno_A$ and $\mu(A_n) \to 1=\mu^*(A)$, but there is no $A 
\in \mathcal{P}(\mathbb{N})$ such that $\mu^*(A) = 1$. This is a 
contradiction, and so the space $L^p$ is not complete, with respect to
the measure $\mu$ defined in (\ref{mumu}).
\end{ex}
%=================
\section{Applications}\label{app}%\mg{app}
\subsection{Inequalities}
In order to investigate some main inequalities in $L^p$ spaces, 
 the following notions and results will be given.
\begin{deff}\label{order-uc} %\mg{order-uc}
\rm
The function $f:\mathbb{R}^+\to\mathbf{X}$ 
is said to be \textit{uniformly continuous}  iff there are a positive element $u\in \mathbf{X}$ and two $(o)$-sequences
$(\sigma_p)_p$ and  $(\delta_p)_p$ in $\mathbf{X}$ and 
$\mathbb{R}^+$, respectively,  such that
$|f(t_1)-f(t_2)|\leq \sigma_p u$ whenever $t_1$, $t_2 \in  
\mathbb{R}^+$  and $p \in {\mathbb N}$ satisfy
$|t_1-t_2| \leq \delta_p$.
\end{deff}
\begin{deff}\label{ud} \rm(\cite[Definition 3.11]{dallas}) %\mg{ud}
Let $I \subset \mathbb{R}$ be a connected set, and
$f:I \to \mathbf{X}$. 
The function $f$ is said to be \textit{uniformly differentiable on $I$} iff there exist
a bounded function $g:I \to \mathbf{X}$ and
two $(o)$-sequences, $(\sigma_p)_p$  and  $(\delta_p)_p$ in
$\mathbf{X}$ and $\mathbb{R}^+$, respectively, with
$$\bigvee \Bigl\{\Bigl|\frac{f(v)-f(u)}{v-u}-g(x)\Bigr|:
u \leq x \leq v , 0 < v-u \leq \delta_p \Bigr\} \leq \sigma_p e$$
for every $p \in \mathbb{N}$. In this case  $g$ is said to be the \textit{uniform derivative} of $f$ in $I$.
\end{deff}
\begin{deff}\label{dconvex} %\mg{dconvex} 
\rm A function 
$f : \mathbf{X}_1 \to \mathbf{X}$ is said to be  a \textit{convex function}, 
iff for every $v \in \mathbf{X}_1$ there exists $\beta_v \in
\mathbf{X}_1$ with 
%\mg{support}
\begin{equation}\label{support}
 f(s) \geq f(v) + \beta_v (s-v)
\end{equation}
for each $s \in \mathbf{X}_1$.
When $\mathbf{X}_1 = \mathbb{R}$, this means 
\[ f(t) \leq f(t_1) + \dfrac{f(t_2) - f(t_1)}{t_2 - t_1} (t-t_1) \quad \mbox{whenever } t, t_1, t_2 \in \mathbb{R},
t_1 < t < t_2.\]
\end{deff}
Let $[a,b]$ be a compact subinterval of the real line.
If $f:[a,b] \to \mathbf{X}$ is uniformly 
continuous on $ [a,b]$, then  the following 
characterization of convex 
vector lattice-valued functions will be obtained.
\begin{thm}\label{semiconv} %\mg{semiconvex}
Assume that $f:[a,b]\to \mathbf{X}$ is a uniformly 
continuous function. Then $f$ is convex if and only if
for every  $a \leq \alpha < \beta \leq b$ it is
%\mg{midpoint} 
\begin{eqnarray}\label{midpoint}
f\Bigl(\dfrac{\alpha+\beta}{2}\Bigr)
\leq \dfrac{f(\alpha)+f(\beta)}{2}.
\end{eqnarray} 
\end{thm} 
\begin{proof}
We begin with proving the "only if'' part.
We can show that (\ref{support}) implies that
$f(c\alpha+(1-c)\beta)\leq cf(\alpha)+(1-c)f(\beta)$
for all $\alpha, \beta$ in $[a,b]$ and $c\in [0,1]$.
Then clearly (\ref{midpoint}) will follow immediately, taking 
$c=1/2$.
To prove the claim, fix $\alpha, \beta, c$ as above, and take 
$v=c\alpha+(1-c)\beta$. Then from (\ref{support}) it follows
$f(\beta)\geq f(v)+\beta_v c (\beta-\alpha),$
and $f(\alpha)\geq f(v)+\beta_v (1-c)(\alpha-\beta),$
by choosing first $s=\beta$ and then $s=\alpha$. 
Now it is
$(1-c)f(\beta)\geq (1-c)f(v)+\beta_v c(1-c)(\beta-\alpha),$
and
$cf(\alpha)\geq c f(v)+\beta_v c (1-c) (\alpha-\beta).$
Summing up the two last inequalities, the implication follows.\\
We now prove the "if'' part.
By usual techniques, from (\ref{midpoint}) it follows that
$f(q\alpha+(1-q)\beta)\leq q f(\alpha)+(1-q)f(\beta)$
for each $\alpha, \beta\in [a,b]$ with $\alpha<\beta$, and for any 
dyadic rational number $q\in [0,1]$.
Since the set of all dyadic rationals of $[0,1]$ is
dense in $[0,1]$, by continuity of $f$ 
it is not difficult to deduce that
%\mg{ineq} 
\begin{eqnarray}\label{ineq}
f(c\alpha+(1-c)\beta)\leq c f(\alpha)+(1-c)f(\beta)
\end{eqnarray}
for any $\alpha, \beta\in [a,b]$ with $\alpha<\beta$, and  
$c\in [0,1]$.
\acr
Now we claim that $f$ is convex. Indeed, if $t$,
$t_1$, $t_2\in [a,b]$ and $t_1<t<t_2$,
one can write $t=\delta t_1+(1-\delta) t_2$, 
where $\delta=
\dfrac{t_2-t}{t_2-t_1}$. A simple application of the inequality in (\ref{ineq})
yields
$ f(t) \leq f(t_1) + \dfrac{f(t_2) - f(t_1)}{t_2 - t_1} (t-t_1).$
This ends the proof.
\end{proof}
Some inequalities 
for convex vector lattice-valued functions will be proven now.
\begin{thm}\label{ucc} %\mg{ucc}
If $f:[a,b] \to \mathbf{X}$ is uniformly 
continuous and $\varphi: \mathbf{X} \to \mathbf{X}$ 
is convex, then $\varphi \circ f$ is uniformly continuous too.
\end{thm}
\begin{proof}
By uniform continuity of $f$ there are a positive element $u\in 
\mathbf{X}$ and two $(o)$-sequences
$(\sigma_p)_p$, $(\delta_p)_p$ in $\mathbf{X}$ and
$\mathbb{R}^+$, respectively, such that
$|f(t_1)-f(t_2)|\leq \sigma_p u$ for all $t_1$, $t_2 \in  \mathbb{R}^+$  and $p \in {\mathbb N}$ satisfying
$|t_1-t_2| < \delta_p$.
Using (\ref{support}) with $v=f(t_1)$, it is
$ \varphi (f(t_2)) \geq \varphi(f(t_1)) + \beta_{f(t_1)} (f(t_2) - f(t_1)),$
while for $v=f(t_2)$ one has
$\varphi(f(t_1)) \geq \varphi (f(t_2)) + \beta_{f(t_2)} (f(t_1) - f(t_2)).$
Setting $\alpha := |\beta_{f(t_1)}| \vee |\beta_{f(t_2)}|$, finally it follows
$
|\varphi(f(t_1)) - \varphi (f(t_2))| \leq \alpha |f(t_1)-f(t_2)|
\leq \alpha \, u \, \sigma_p 
$
whenever $|t_1-t_2| < \delta_p$.
\end{proof}
\begin{thm} { \rm (Jensen inequality)}
If $f:[a,b] \to \mathbf{X}$ is uniformly continuous and 
$\varphi: \mathbf{X} \to \mathbf{X}$ is convex, then
\[  \varphi \Bigl( \int_a^b f(t) d\mu(t) \Bigr) 
\leq \int_a^b \varphi(f(t)) d\mu(t)\]
whenever $\mu([a,b]) = 1$.
\end{thm}
\begin{proof}
By \cite[Proposition 3.12]{dallas} and  Proposition \ref{ucc}, both $f$ and $\varphi \circ f$ 
are bounded and integrable.
Let $\displaystyle{
\tau:= \int_a^b f(t) d\mu(t)}$ and $m, M \in \mathbf{X}$ 
be such that $m \leq f(t) \leq M$ for every $t \in [a,b]$.
Let $\beta_{\tau}$ be the element associated to $\tau$ 
according to Definition \ref{dconvex}. Then
\begin{eqnarray*}
\varphi(f(t)) \geq \varphi(\tau) + \beta_{\tau} (f(t) - \tau).
\end{eqnarray*}
Integrating both members, it is
\begin{eqnarray*}
\int_a^b \varphi(f(t)) d\mu(t) &\geq& \varphi \Bigl(\int_a^b f(t) d\mu(t) \Bigr) +
\int_a^b\beta_{\tau} (f(t) - \tau) d\mu(t) =\\
&=& \varphi \Bigl (\int_a^b f(t) d\mu(t) \Bigr),
\end {eqnarray*}
that is the assertion.
\end{proof}
\begin{thm}\label{f1up} %\mg{f1up}
Let $f : [a,b] \to \mathbf{X}$ be a uniformly differentiable and 
convex function, and set
%\mg{r} 
\begin{eqnarray}\label{r}
r(t):=f\Bigl( \dfrac{a+b}{2} \Bigr) + f^{\prime} 
\Bigl(\dfrac{a+b}{2}\Bigr) \Bigl(t - \dfrac{a+b}{2}\Bigr),
\quad t \in [a,b].
\end{eqnarray}
Then $f(t) \geq r(t)$ for every $t\in [a,b]$.
\end{thm}
\begin{proof}
Fix arbitrarily $t_1$, $t_2 \in [a,b]$ with $t_1 < t_2$.
By convexity of $f$, for every $t \in ]t_1, t_2[$, it is
%\mg{fine} 
\begin{eqnarray}\label{fine}
&& \dfrac{f(t)-f(t_1)}{t-t_1} \leq  \dfrac{f(t_2)-f(t_1)}{t_2 -
t_1}, \quad \dfrac{f(t)-f(t_2)}{t-t_2} 
\geq  \dfrac{f(t_2)-f(t_1)}{t_2 -t_1}. 
\end{eqnarray}
Thanks to uniform differentiability of $f$, 
(\ref{fine}) implies
%\mg{fine2} 
\begin{eqnarray}\label{fine2}
f^{\prime} (t_1) \leq  \dfrac{f(t_2)-f(t_1)}{t_2 -t_1} \leq 
f^{\prime}(t_2).
\end{eqnarray}
Using the first inequality in (\ref{fine2}) with 
$t_1 = \dfrac{a+b}{2}$, $t_2=t$, it follows
\[ f(t) \geq  f\Bigl( \dfrac{a+b}{2} \Bigr) + f^{\prime} 
\Bigl(\dfrac{a+b}{2}\Bigr) \Bigl(t - \dfrac{a+b}{2}\Bigr)= r(t).\]
By applying the second inequality in (\ref{fine2}) with
$t_1=t$, $t_2 = \dfrac{a+b}{2}$, we obtain again
$f(t) \geq r(t)$. The assertion follows from arbitrariness
of $t_1$ and $t_2$. 
\end{proof}
\begin{thm}\label{HH} %\mg{HH} 
{\rm (Hermite-Hadamard inequality)}
Let $f : [a,b] \to \mathbf{X}$ be a uniformly differentiable and convex function. Then
\[ f\Bigl( \dfrac{a+b}{2} \Bigr) \leq \dfrac{1}{b-a} \int_a^b f(t) d\mu(t) \leq \dfrac{f(a)+f(b)}{2}.\]
\end{thm}
\begin{proof}
Let $r$ be as in (\ref{r}), and set
%\mg{s}
 \begin{eqnarray}\label{s}
r^*(t) = f(a) + \dfrac{f(b) - f(a)}{b-a} (t-a),
\quad t \in [a,b].
\end{eqnarray}
By \cite[Proposition 3.12]{dallas}, $f \in L^1(\lambda)$ and 
\begin{eqnarray*}
\int_a^b r(t) d\mu(t) = (b-a) f\Bigl( \dfrac{a+b}{2} \Bigr) \leq 
\int_a^b f(t) d\mu(t).
\end{eqnarray*}
Since $f$ is convex,
$f(t) \leq r^*(t)$ for every $t\in [a,b]$
and
\[  \int_a^b f(t) d\mu(t) \leq \int_a^b r^*(t) d\mu(t) = (b-a)  \dfrac{f(a)+f(b)}{2}.\]
\end{proof}
\begin{thm}\label{Fejer} %\mg {Fejer} 
{\rm (F\'{e}jer  inequality)}
Let $f : [a,b] \to \mathbf{X}$ be a uniformly differentiable and convex function and let
 $w:[a,b] \to \mathbb{R}_0^+$
be a uniformly continuous map such that $w(a+t) = w(b-t)$ for every $0 \leq t \leq \dfrac{a+b}{2}$. Then 
 \[ f\Bigl( \dfrac{a+b}{2} \Bigr) \int_a^b w(t) d\mu(t) \leq 
 \int_a^b f(t)w(t) d\mu(t) \leq \dfrac{f(a)+f(b)}{2} 
\int_a^b w(t) d\mu(t).\]
\end{thm}
\begin{proof}
Again by \cite[Proposition 3.12]{dallas},
$fw$ is in $L^1(\lambda)$. Let $r$, $r^*$ be as in 
(\ref{r}) and (\ref{s}), respectively. Then 
%\mg{rsw} 
\begin{eqnarray}\label{rsw}
r(t)w(t) \leq f(t)w(t) \leq r^*(t)w(t), \quad t \in [a,b].
\end{eqnarray}
By integrating all members in (\ref{rsw}) 
the assertion follows, since 
$$\dfrac{b+a}{2} \int_a^b w(t) 
d\mu(t) = \int_a^b t w(t) d\mu(t).$$
\end{proof}

In this setting the case $p=2$ is particularly interesting. Indeed it is: 
\begin{thm}\label{s-ineq} %\mg{s-ineq} 
{\rm (Schwartz inequality)}
If $f$, $g \in L^2$, then $fg$ is integrable. Moreover it is
\begin{eqnarray*}
&& \Bigl( \int_G |f(t)g(t)| d\mu(t) \Bigr)^2 \leq \Bigl( \int_G (f(t))^2 
d\mu(t) \Bigr) \Bigl( \int_G (g(t))^2 d\mu(t) \Bigr),\\
&&
\sqrt{\Bigl\| \Bigl( \int_G |f(t)g(t)| d\mu(t) \Bigr)^2 \Bigr\|_e }
\leq \| f\|_2 \| g \|_2.
\end{eqnarray*}
\end{thm}
\begin{proof}
Since $f$, $g \in L^2$, there exist defining sequences $(f_n)_n$ and $(g_n)_n$, related to $f$ and $g$, such that $(f^2_n)_n$ and $(g^2_n)_n$ are defining for $f^2, g^2$ respectively.
Then, $f+g \in L^2$ and $(f_n + g_n)^2$ is defining for $(f+g)^2$ (see also Proposition \ref{sum-p}). This means that 
\( f_n g_n := 1/2 ( (f_n + g_n)^2 - f^2_n - g^2_n ) \), $n\in 
\mathbb{N}$, is a defining sequence for $fg$ and this shows that $fg$ 
is integrable. Observe also that all involved
mappings can be taken non negative. 
Now, in order to prove the Schwartz inequality, thanks to the definition of integral, it is sufficient to prove it just for simple
 functions $f$ and $g$.
So, assume that $$\displaystyle{f= \sum_{i=1}^n c_i 1_{E_i}}, \quad
\displaystyle{g = \sum_{i=1}^n d_i 1_{E_i}}$$ (without loss of generality, we can take
the same partition for both mappings). Now
\begin{eqnarray*}
&& \int_G f(t)g(t) d\mu(t) 
= \sum_{i=1}^n c_i d_i \mu(E_i), \quad
\int_G (f(t))^2 d\mu(t) 
= \sum_{i=1}^n c^2_i  \mu(E_i), \\
&& \int_G (g(t))^2 d\mu(t) = \sum_{i=1}^n  d^2_i \mu(E_i).
\end{eqnarray*}
Therefore, 
\begin{eqnarray*}
\Bigl( \int_G f(t)g(t) d\mu(t) \Bigr)^2 =
\sum_{i=1}^n c^2_i d^2_i \mu(E_i) + 2 \sum_{i < j} c_i c_j d_i d_j \mu(E_i) \mu(E_j),
\end{eqnarray*}
\begin{eqnarray*}
&& \Bigl( \int_G (f(t))^2 d\mu(t) \Bigr)^2 
\Bigl( \int_G (g(t))^2 d\mu(t) \Bigr)^2 =
\sum_i \sum_j c^2_i d^2_j \mu(E_i) \mu(E_j) +
\\ &&+ \sum_{i < j} c^2_i d^2_j\mu(E_i) \mu(E_j) +
\sum_{i < j} c^2_j d^2_i \mu(E_i) \mu(E_j).
\end{eqnarray*}
Comparing the last two formulas, it is clear that
\begin{eqnarray*}
&& \Bigl( \int_G (f(t))^2 d\mu(t)  \Bigr)^2 \cdot
\Bigl( \int_G (g(t))^2 d\mu(t) \Bigr)^2 
- \Bigl( \int_G f(t)g(t) d\mu(t) \Bigr)^2 = \\ &&=
\sum_{i < j} (c^2_i d^2_j + c^2_j d^2_i - 2 c_i c_j d_i d_j)\mu(E_i) \mu(E_j).
\end{eqnarray*}
Since $c^2_i d^2_j + c^2_j d^2_i - 2 c_i c_j d_i d_j= 
(c_i d_j - c_j d_i)^2$ for all $i$, $j$, easily it follows
\[
\Bigl( \int_G f(t)g(t) d\mu(t) \Bigr)^2  
\leq \Bigl( \int_G (f(t))^2 d\mu(t) \Bigr)^2 \Bigl( 
\int_G (g(t))^2 d\mu(t) \Bigr)^2.
\]
Taking the norm, we finally deduce
\begin{eqnarray*}
\Bigl\| \Bigl( \int_G f(t)g(t) d\mu(t) 
\Bigr)^2 \Bigr\|_e  &\leq& 
\Bigl\| \Bigl( \int_G (f(t)^2 d\mu(t)  \Bigr)^2 \Bigl( \int_G 
(g(t))^2 d\mu(t) \Bigr)^2 \Bigr\|_e \leq \\&\leq&
\Bigl\| \Bigl( \int_G (f(t))^2 d\mu(t) \Bigr)^2 
\Bigr\|_e  \cdot \Bigl\| \Bigl( \int_G 
(g(t))^2 d\mu(t) \Bigr)^2 \Bigr\|_e,
\end{eqnarray*}
where
the last inequality follows from the idempotence of $e$. In conclusion
\[ \sqrt{\Bigl\| \Bigl( 
\int_G f(t)g(t) d\mu(t) \Bigr)^2 \Bigr\|_e } \leq \| f\|_2 \|g \|_2.
\]
\end{proof}

\subsection{The Brownian Motion}
In order to obtain a concrete application in Stochastic Integration, we 
assume that $\displaystyle{\mathbf{B}:=(B_t)_{0\leq
t<T}:  [0,T] \to L^2(\Omega)}$ (with $T<+\infty$) is the standard 
Brownian Motion defined on a probability space $(\Omega,\Sigma,P)$. 
As observed in \cite{dallas},  $L^2$ has not an order unit, 
but thanks to the well-known Maximum Principle (see e.g. \cite{brei})  
there is a positive element $Z\in L^2$ with $|B(t)|\leq Z$ for all $t\in [0,T]$. 
Moreover,  there exists a positive random variable $W$ in $L^2$ such that 
$|B(t+h)-B(t)|\leq |h|^{1/4} W$
whenever $t,t+h\in [0,T]$  (see e.g. \cite{garsia, nualart}). 
Thus, taking  $\mathbf{X}$ as the (complete) subspace of $L^2$ 
generated by all elements dominated by some real multiple of $W+Z$, we see that $\mathbf{X}$
 has an order unit (i.e. $W+Z$), and $\textbf{B}$ is a uniformly continuous 
 $\mathbf{X}$-valued function defined on 
$[0,T]$. \acr

\begin{ex}\label{4.12} %\mg{4.12} 
\rm  
Fix now $T> a > 1$, and set
$$f(t):=\Bigl\{\begin{array}{ll} 
0, & \text{if  } 0\leq t\leq a \quad {\rm or} \quad  t\geq T, \\
(t-T)(B_t-B_a), & \text{if  } a\leq t\leq T.
\end{array}
$$
Note that $f$ is a uniformly continuous 
function from $[0,+\infty[$ to $L^2(\Omega)$. Then 
%\mg{Fi} 
\begin{eqnarray}\label{Fi}
\Phi(x,s)=\frac{x}{s^x}\int_0^sf(t)t^{x-1}d\mu(t)
\end{eqnarray}
is a solution of the partial differential equation
$s\dfrac{\partial^2 \Phi}{\partial x \partial s}+x
\dfrac{\partial \Phi}{\partial x}+\Phi(x,s)=f(s).$
Indeed, it is sufficient to differentiate with respect to $x$ the equation 
\begin{eqnarray}{\label{riferimento}}
s\frac{\partial{\varphi}(x,s)}{\partial s}+x\varphi(x,s)=x f(s),
\end{eqnarray}
taking into account that it admits the solution
$\varphi(x,s)=\dfrac{1}{s^{x}}\Bigl(c+x\int_0^sf(t)t^{x-1}d\mu(t)\Bigr),$
where $c$ is an arbitrary real constant. \acr

A more general solution for the previous 
partial differential equation is
$$\Phi(x,s)=F(s)+\frac{G(x)}{s^x}+\frac{x}{s^x}\int_0^sf(t)t^{x-1}d\mu(t),$$
where $F$ and $G$ are arbitrary sufficiently regular functions.
\vspace{2mm}

In the following two figures,  the graphs of the square of the Brownian 
Bridge $f$ and its corresponding surface $\varphi(x,s)$
are represented, respectively.
\begin{figure}[htbp]
\centering
\fbox{\includegraphics[width=.48\textwidth]{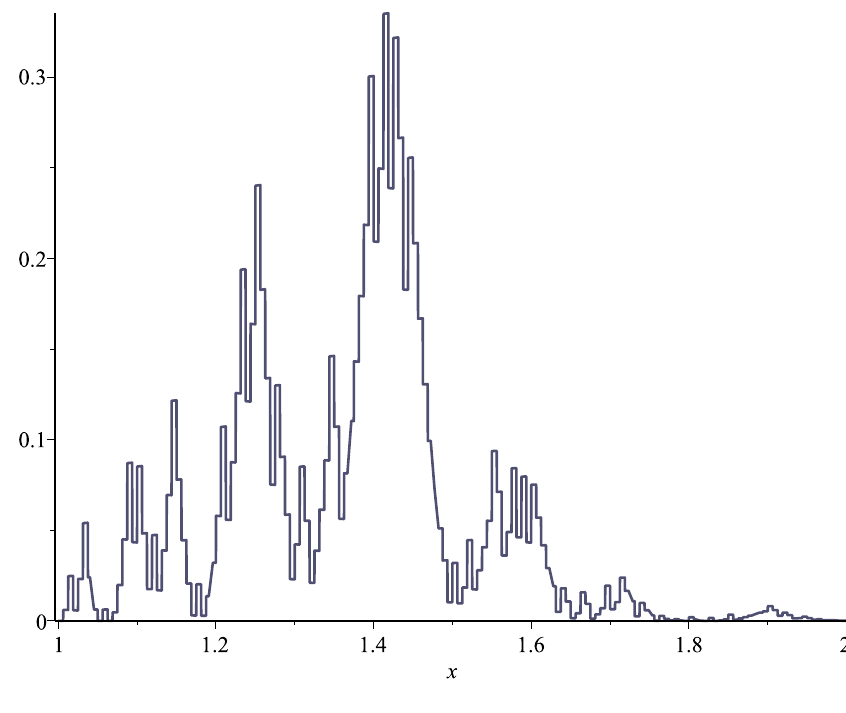}}
\caption{The $f^2$ graph.}\quad
\end{figure}
\begin{figure}[htbp]
\centering
\fbox{\includegraphics[width=.48\textwidth]{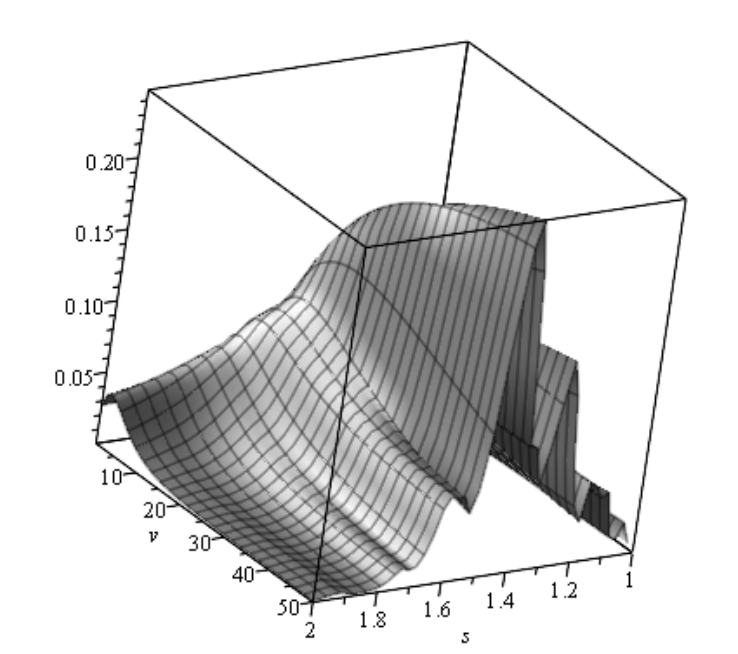}}
\caption{The surface $\varphi(x,s)$ associated with $f^2$.}
\end{figure}
\end{ex}

The two graphs show also that the surface $\varphi(x,s)$
uniformly converges to $f(s)$ as $x$ tends to $+\infty$, according to
\cite[Proposition 3.13]{dallas}: this could be verified, comparing the 
graph of Figure 1 with the section of the surface in Figure 2 
corresponding to $x=50$. \acr

The previous example suggests a different type of operator, 
acting with respect to stochastic differentials.
Indeed, with the same notations as above,
\begin{ex}\rm
Let 
$$\displaystyle{
\Psi(x,s):=\frac{x}{s^x}\int_0^s t^{x-1}df(t)},$$
where now $df(t)$ is the stochastic differential of the {\em Brownian Bridge} $f(t)$, namely
$df(t)=(B_t-B_a)dt+(t-T)dB_t,$
when $a\leq t\leq T$. \acr
The existence of the above integral is ensured by the fact that $t\mapsto t^{x-1}$ is a 
function of class $C^1$, 
and the trajectories of $f(t)$ are continuous 
maps.   Moreover, the well-known formula of integration by parts leads to
\begin{eqnarray*}
\Psi(x,s)=\frac{x}{s^x}\Bigl(s^{x-1}f(s)-\int_0^sf(t)d(t^{x-1})\Bigr)=%\\ &=&
\frac{x}{s}f(s)-\frac{x(x-1)}{s^x}\int_0^sf(t)t^{x-2}dt.
\end{eqnarray*}
So, in conclusion 
%\mg{PhiPsi} 
\begin{eqnarray}\label{PhiPsi}
\Psi(x,s)=\frac{x}{s}(f(s)-\Phi(x-1,s)),
\end{eqnarray}
where $\Phi$ is as in (\ref{Fi}).
This leads to a stochastic differential equation satisfied by $\Psi$. \acr
However, in order to avoid confusion, we  point out that expressions like $d\Psi$ or $d\Phi$ always refer 
to the variable $s$, since $x$ is just a parameter.
\\
Differentiating (\ref{PhiPsi}), it follows
$$d\Psi(x,s)=-\frac{x}{s^2}(f(s)-\Phi(x-1,s))ds+\frac{x}{s}(df(s)-\Phi'(x-1,s)ds).$$
(Here the {\em stochastic differential} involved is $df$, rather than $dB$).
\\
Thanks to (\ref{PhiPsi}) and (\ref{riferimento}), it is
$$d\Psi(x,s)=-\frac{1}{s}\Psi(x,s)ds+\frac{x}{s}df(s)-\frac{x}{s}\frac{(x-1)(f(s)-\Phi(x-1,s))}{s}ds,$$
and, again by (\ref{PhiPsi}),
$\displaystyle{
d\Psi(x,s)=-\frac{1}{s}\Psi(x,s)ds+\frac{x}{s}df(s)-\frac{x-1}{s}\Psi(x,s)ds}$.
In conclusion, after a final simplification, $\Psi(x,s)$ satisfies the stochastic equation
%\mg{stochequa} 
\begin{eqnarray}\label{stochequa}
d\Psi(x,s)=-\frac{x}{s}\Psi(x,s)ds+\frac{x}{s}df(s).
\end{eqnarray}
\acr

\noindent A more general solution of the last equation is
$$\displaystyle{\widetilde{\Psi}(x,s)=\frac{1}{s^x}\Bigl(k+x\int_0^st^{x-1}df \Bigr)},$$
where $k$ is an arbitrary constant.
\end{ex}

Another interesting application of this type of operator can be found in 
detecting a regular signal, when it is distorted by a random noise.
More precisely
\begin{ex}\rm
Let $h:[0,+\infty[\to \erre$ be any $C^1$ map, satisfying the condition 
that $h(x)=h'(x)=0$ for all $x\in [0,a]$ and $x\geq T$, where $a,T$ are 
fixed positive numbers, $1<a<T$.
 As above, let $f$ denote the 
Brownian Bridge, like in  Example \ref{4.12}, and fix any positive 
number $\varepsilon$. Then set $G(s)=h(s)+\varepsilon f(s)$, for $s\in 
[0,+\infty[$, and 
$$
\displaystyle{\Psi_G(x,s)=\frac{x}{s^x}\int_0^s t^{x-1}dG(s)},$$
i.e. 
$$\Psi_G(x,s)=\frac{x}{s^x}\int_0^s t^{x-1}h'(t)dt+
\frac{x}{s^x}\int_0^s \varepsilon t^{x-1}df(t).$$
Thanks to (\ref{riferimento}) and (\ref{PhiPsi}), 
$$\Psi_G(x,s)=\Phi_{h'}(x,s)+\frac{\varepsilon x}{s}\Bigl( f(s)-\Phi_f(x-1,s) \Bigr),$$
where the operator $\Phi$ has the same meaning as in Example \ref{4.12}.
Now, thanks to \cite[Proposition 3.13]{dallas}, 
$\displaystyle{
\lim_{x\to +\infty}\Phi_{h'}(x,s)=h'(s)}$
and 
$$\displaystyle{\lim_{x\to +\infty} (f(s)-\Phi_f(x-1,s))=0},$$
both uniformly with respect to $s$. So, choosing 
$x=\varepsilon^{-1}$ in the above quantities, it is
$$\displaystyle{
\lim_{\varepsilon\to 0}\Psi_G (\varepsilon^{-1},s )=h'(s)}.$$
\end{ex}
%=======================
\section*{Conclusions}
After an abstract introduction of $L^p$ spaces for vector 
lattice-valued functions, some applications are found, first
extending some classical inequalities to the vector lattice setting, and then finding approximations
 for stochastic processes like
the Brownian Motion and the Brownian Bridge, based on the moment operator; 
as an outcome,  
this method leads to solve some stochastic differential equations.
%====================== BIBLIOGRAFIA=========================

\end{document}